\documentclass[12pt]{article}
\usepackage{setspace}
\usepackage{pstricks}
\usepackage{graphicx}
\usepackage[left=3cm,top=3cm,bottom=3cm,right=3cm]{geometry}
\usepackage{amsthm}
\usepackage{amssymb}
\theoremstyle{definition}

\newtheorem*{prf}{Proof}

\newtheorem{thm}{Theorem}
\newtheorem{defn}{Definition}
\newtheorem{lem}{Lemma}

\title{Countable chains and infinite joins in effectively closed sets of Cantor space}
\author{Ahmet \c{C}evik\footnote{E-mail: a.cevik@hotmail.com}\\ {\small Gendarmerie and Coast Guard Military Academy, Ankara, Turkey}}
\date{ }
\begin{document}
\maketitle

\begin{abstract}
We prove that there exists a countable infinite sequence of non-empty special $\Pi^0_1$ classes $\{\mathcal{P}_i\}_{i\in\omega}$ such that no infinite union of elements of any $\mathcal{P}_i$ computes the halting set. We then give a generalized form of lower and upper cone avoidance for infinite unions. That is, we show that for any special $\Pi^0_1$ class $\mathcal{P}$ and any countable sequence of sets in $\mathcal{P}$, $\mathcal{P}$ has a member that is not computable by the infinite union of elements of the sequence. We also prove the upper cone counterpart, that for any non-recursive set $X$, every non-empty $\Pi^0_1$ class contains a countable sequence of members whose join does not compute $X$. We finally show that there exists a $\Pi^0_1$ class whose degree specrum is a countably infinite strict chain.
\end{abstract}



\section{Introduction}

Effectively closed sets, i.e.,  $\Pi^0_1$ classes, has been a central theme in classical recursion theory. More specifically, the problem of determining degree theoretic complexity of members of $\Pi^0_1$ classes, going back to Kleene \cite{SK}, has resulted in a rich theory. Since the Cantor space is compact, degree theoretic complexity of members of $\Pi^0_1$ classes also determines what kind of reals can be defined by compactness rather than using replacement. One seminal paper on this subject is by Jockush and Soare \cite{JS1} \cite{JS2}, who showed some very interesting degree theoretic properties of members of $\Pi^0_1$ classes. Many of these results are known as {\em basis theorems} for $\Pi^0_1$ classes. A typical basis theorem tells us that every $\Pi^0_1$ class has a member, or a member of degree, of a particular kind. It may also be the case that not every $\Pi^0_1$ class has members with the desired property. We prove that there exists a countable infinite sequence of non-empty special $\Pi^0_1$ classes $\{\mathcal{P}_i\}_{i\in\omega}$ such that no infinite union of elements of any $\mathcal{P}_i$ computes the halting set. We then give a generalized form of lower and upper cone avoidance for infinite unions. That is, we show that for any special $\Pi^0_1$ class $\mathcal{P}$ and any countable sequence of sets in $\mathcal{P}$, $\mathcal{P}$ has a member that is not computable by the infinite union of elements of the sequence. We also give the upper cone counterpart, that for any non-recursive set $X$, every non-empty $\Pi^0_1$ class contains a countable sequence of members whose join does not compute $X$. We finally show that there exists a $\Pi^0_1$ class whose degree specrum is a countably infinite strict chain.

\subsection{Notation}
We let $\omega$ denote the set of natural numbers. Let $2^{<\omega}$ denote the set of all finite sequences over $\{0,1\}$. We denote sets of natural numbers with capital Latin letters such as $A,B,C,D$. The subset relation is denoted by $\subset$. We identify a set $A\subset\omega$ with its characteristic function $f:\omega\rightarrow\{0,1\}$ such that, for any $n\in\omega$, if $n\in A$ then $f(n)=1$, and if $n\not\in A$ then $f(n)=0$.
\vspace{0.5cm}

{\bf Turing functionals.}
We can have $\{\Psi_i\}_{i\in\omega}$ as an effective enumeration of the Turing functionals. $\Psi_e$ is {\em total} if it is defined for every argument, otherwise it is called {\em partial}. For any $A\subset\omega$ and $n\in\omega$, $\Psi_e(A;n)\downarrow=m$ denotes that the $e$-th Turing functional with oracle $A$ on argument $n$ is defined and equal to $m$. For any $A$ and $n$, $\Psi_e(A;n)\uparrow$ denotes it is not the case that $\Psi_e(A;n)\downarrow$. Since $\Psi_e(A)$ denotes a partial function and since we identify subsets of $\omega$ with their characteristic functions, it is reasonable to write $\Psi_e(A)=B$ for some $B\subset\omega$. If a set of natural numbers $A$ is computable in $B$, we denote this by $A\leq_T B$. If neither $A\leq_T B$ nor $B\leq_T A$, we say that $A$ and $B$ are {\em Turing incomparable}. If $A\leq_T B$ and $B\leq_T A$, then $A$ and $B$ are {\em Turing equivalent}, i.e., $A\equiv_T B$. We denote Turing degrees with boldcase letters ${\bf a,b,c}$. 
Partial functions are also denoted by $f,g$. We let $\left\langle i,j\right\rangle$ be a computable bijection $\omega\times\omega\rightarrow\omega$. 
\vspace{0.5cm}

{\bf Strings.}
We denote finite strings with lowercase Greek letters like $\sigma,\tau,\eta,\rho,\pi,\upsilon$. We let $\sigma*\tau$ denote the concatenation of $\sigma$ followed by $\tau$. We let $\sigma\subset\tau$ denote that $\sigma$ is an initial segment of $\tau$. We say a string $\sigma$ is {\em incompatible with $\tau$} if neither $\sigma\subset\tau$ nor $\tau\subset\sigma$. Otherwise we say that $\sigma$ is {\em compatible} with $\tau$. Similarly, we say that $\sigma$ is an {\em extension} of $\tau$ if $\tau\subset\sigma$. Let $|\sigma|$ denote the length of $\sigma$. We let $\sigma(i)$ denote the $(i+1)$st bit of $\sigma$. 
\vspace{0.5cm}

{\bf Use function conventions.}
For any $\sigma\in 2^{<\omega}$ and for any $n\in\omega$, we let $\Psi_e(\sigma;n)$ be defined and equal to $\Psi_e(A;n)$ if $\sigma(i)=A(i)$ for all $i<|\sigma|$ and if computing $\Psi_e(A;n)$ requires only values $A(i)$ for $i<|\sigma|$. 
$\Psi_i(\sigma)[s]$ denotes $\Psi_i(\sigma)$ defined at stage $s$. For convenience, we suppose that $\Psi_i(\sigma)[s]$ is defined only on numbers $n\leq s$. For a set $A\subset\omega$, we define the {\em jump} of $A$, denoted $A'$, to be $\{e:\Psi_e(A;e)\downarrow\}$. The {\em halting set} is $K=\{e:\Psi_e(e)\downarrow\}$. 
\vspace{0.5cm}

{\bf Trees and $\Pi^0_1$ classes.}
A set $\Lambda$ of strings is {\em downward closed} if $\sigma\in \Lambda$ and $\tau\subset\sigma$ then $\tau\in \Lambda$. Occasionally we refer to downward closed sets of strings as {\em trees}. We say that a set $A$ {\em lies} on $\Lambda$ if there exist infinitely many $\sigma$ in $\Lambda$ such that $\sigma\subset A$. A set $A$ is a {\em path} on $\Lambda$ if $A$ lies on $\Lambda$. So if $A$ is a path on $\Lambda$, then every initial segment of $A$ is in $\Lambda$. We denote the set of infinite paths of $\Lambda$ by $[\Lambda]$. 
We say that a string $\sigma\in T$ is {\em infinitely extendible} if there exists some $A\supset\sigma$ such that $A\in[T]$. 
If $\sigma,\tau\in T$ and $\sigma\subset\tau$ and there does not exist $\sigma'$ with $\sigma\subset\sigma'\subset\tau$ then we say that $\tau$ is an {\em immediate successor} of $\sigma$ in $T$ and $\sigma$ is the {\em immediate predecessor} of $\tau$ in $T$. 

We say that $\mathcal{P}\subset 2^\omega$ is a {\em $\Pi^0_1$ class} if there exists a downward closed computable set of strings $\Lambda$ such that $\mathcal{P}=[\Lambda]$. We can then have an effective enumeration $\{\Lambda_i\}_{i\in\omega}$ of downward closed computable sets of strings such that for any $\Pi^0_1$ class $\mathcal{P}$ there exists some $i\in\omega$ such that $\mathcal{P}$ is the set of all infinite paths through $\Lambda_i$. Call a $\Pi^0_1$ class {\em special} if it does not contain a computable member. The {\em degree spectrum} of a $\Pi^0_1$ class $\mathcal{P}$ is defined by the set of all degrees ${\bf a}$ such that there exists a member $A\in\mathcal{P}$ of degree ${\bf a}$. Occasionally we may use subsets of Baire space $\omega^\omega$ but unless we explicitly state that, we will be working in Cantor space. 

It is important to note that the compactness property of the Cantor space is provided by the Weak K\"{o}nig's Lemma which tells us that if $\Lambda$ is an infinite downward closed set of finite strings, then there exists an infinite path through $\Lambda$.

For a detailed account of computability, the reader may refer to \cite{Cooper} or \cite{Soare}. An extensive survey for $\Pi^0_1$ classes can be found in \cite{CenzerRec}, \cite{Diamondstone} and \cite{Downey}.

\section{Finite joins of members of infinitely many classes}

We construct an infinite sequence of effectively closed sets in which no finite join of members of any of these sets computes the halting set. We need to generalize the join operator to any finite number of sets.

\begin{defn}
Let $A_0,A_1,\ldots A_n$ be a sequence of sets of natural numbers. Then, the {\em finite join} is defined as $\bigoplus_{i\leq n} A_i=A_0\oplus A_1\oplus \ldots \oplus A_n=\bigcup_{j\leq n}\{(n+1)i+j : i\in A_j\}$.
\end{defn}

The theorem we prove now will show that we can define infinitely many $\Pi^0_1$ classes from which if we choose finitely many members from any of them, the join of these members still does not compute the halting set. 

\begin{thm}
There exist a countably infinite sequence of non-empty special $\Pi^0_1$ classes $\{\mathcal{P}_n\}_{n\in\omega}$ such that for any finite sequence $A_0,A_1,\ldots A_n$ of elements, each of which is a member of some $\mathcal{P}_k$, $\bigoplus_{i\leq n} A_i$ does not compute the halting set.
\end{thm}

\begin{prf}
We want to satisfy $\emptyset'\not\leq_T \bigoplus A_i$ by defining a countable sequence of sets $\{T_i\}_{i\in\omega}$ via recursive approximation such that $\mathcal{P}_i=[T_i]$, and a set $D$ such that $D\not\leq_T \bigoplus A_i$. We define $T_i$'s using dovetailing, so for convenience we may assume that $T_i$ is not defined until stage $i$.

The requirements are as follows:
\vspace{0.5cm}

$R_{\langle {2e+1},i\rangle}$: If $S\in\mathcal{P}_i$, then $S\neq\Psi_e$

$R_{2e+2}$: If $A_i\in\mathcal{P}_j$ for some $j\in\omega$, then $\Psi_e(\bigoplus A_i)\neq D$.
\vspace{0.5cm}

Construction.

At stage $s=0$, we define $\emptyset$ to be in every $T_i$.

At stage $s>0$
\begin{enumerate}
\item[(i)] For every $i\leq s$ find the least string, if exists, $\tau\in T_i$, such that $\tau$ is of level $2e+1$ and $\tau\subset\Psi_e[s]$. Let $\tau_0\in T_i$ be the immediate predecessor of $\tau$ and let $\tau_1$ be a leaf of $T_i$ extending $\tau_0$ and incompatible with $\tau$. Stop enumerating any strings extending $\tau_0$ in $T_i$, i.e. declare them terminal.

\item[(ii)] If the enumerated strings are at even level, we consider all $n$-tuples of mutually incompatible strings $\langle \sigma_0,\sigma_1,\ldots,\sigma_n\rangle$ of that level, each of which is in some $T_i$, that we have not yet declared to be terminal. If not done already, we also assign them a unique follower $x$ not yet enumerated in $D$. Then, at any subsequent stage, if we see some $\sigma_k'\supset\sigma_k$ for all $k\leq n$ already in some $T_j$, $j\leq s$, at level $2e+2$ that are not yet declared to be terminal such that $\Psi_e(\sigma_0'\oplus\sigma_1'\oplus\ldots\sigma_n';x)\downarrow=0$ and $D(x)=0$, we enumerate $n$ into $D$. Then we declare all other extensions of $\sigma_0,\sigma_1,\ldots,\sigma_n$ in $T_i$, other than $\sigma_0',\sigma_1',\ldots,\sigma_n'$, to be terminal.
\end{enumerate}

After these instructions, choose two incompatible strings extending each leaf of $T_j$ for all $j\leq s$ that have not yet declared to be terminal, and enumerate them into $T_j$.
\vspace{0.5cm}

\subsection{Verification}

We shall first show that for each $i\in\omega$, $[T_i]$ is a non-empty $\Pi^0_1$ class. In the construction, $T_i$ may not be necessarily downward closed. We need to show explicitly that there exists a downward closed computable set of strings $\Lambda_i$ such that $[T_i]=[\Lambda_i]$. We let $\Lambda_i$ be the set of all strings which are initial segments of strings in $T_i$ at any stage $s$. That is, we define $\Lambda_{i,s}$ to be the set of all initial segments of $\sigma$, for any $\sigma\in T_{i,s}$. Let $\Lambda_i=\bigcup \Lambda_{i,s}$. We next show that $\Lambda_i$ is downward closed, computable and $[\Lambda_i]=[T_i]$. Now $\Lambda_i$ is computable since we may assume without loss of generality that $\sigma\in\Lambda_i$ iff $\sigma\in\Lambda_{i,s}$ for $s=|\sigma|$. Clearly, every infinitely extendible string in $T_i$ is also in $\Lambda_i$ by the defintion of $\Lambda_i$. The opposite direction is also true. Suppose that $\sigma$ is not infinitely extendible in $\Lambda_i$. Then $\sigma$ must be a leaf of $T_i$ in which case $\sigma$ is not infinitely extendible in $T_i$ since otherwise $\sigma$ would be infinitely extendible in $\Lambda_i$. 

The non-emptiness of $[T_i]$ follows from the fact that we maintain, at all stages, that if $\sigma$ is not a leaf and there is at least one leaf extending $\sigma$ that has not yet been declared terminal, then there is a pair of incompatible leaves extending $\sigma$ which have not yet been declared terminal. By virtue of construction, at stage $s=0$, two incompatible extensions of $\emptyset$ will be enumerated into $T_i$. At each next stage $s>0$, we enumerate into $T_{i,j}$, for every $j\leq s$, two incompatible strings extending each leaf of $T_{i,j}$. Now since only one of every pair of incompatible strings satisfies the $R_{2e+1}$ or the $R_{2e+2}$ requirement, we guarantee to leave at least one node alive.

\begin{lem}
$R_{2e+1}$ is satisfied.
\end{lem}
\begin{prf}
Given an index $i\in\omega$, suppose that $S\in\left[T_i\right]$ and $S=\Psi_e(\emptyset)$ for some $e$. Then for all $\sigma\subset S$ , where $\sigma\subset\Psi_e(\emptyset)$, we have $\sigma\in T_i$ at level $2e+1$. Let $\sigma_0$ be the immediate predecessor of $\sigma$. Then any extensions of $\sigma_0$, compatible with $\sigma$, are not enumerated into $T_i$. But then, $\sigma\subset\Psi_e(\emptyset)$ is in $T_i$ for finitely many $\sigma$'s. A contradiction.
\end{prf}

\begin{lem}
$R_{2e+2}$ is satisfied.
\end{lem}
\begin{prf}
Suppose the contrary that there exists a finite sequence of sets $A_0,A_1,\ldots A_n$ such that, for $i\leq n$, $A_i\in \mathcal{P}_j$ for some $j\in\omega$ and that $\Psi_e(\oplus A_i)=D$ for some $e$. Then for each $A_i$ there exists some $\sigma_i$ such that $\sigma_i\subset A_i$ in $T_i$ at level $2e$ and there exists $\sigma_i'\supset\sigma_i$ for each $i$ such that $\Psi_e(\oplus\sigma_i';m)=D(m)$, where $m$ is the assigned follower for the tuple $\langle\sigma_0,\sigma_1,\ldots,\sigma_n\rangle$. But then, by virtue of the construction, this is a contradiction. This completes the proof of the theorem.
\end{prf}
\qed
\end{prf}

\section{Cone avoidance up to infinite joins}
In this part we observe that the lower cone and upper cone avoidance theorems given in \cite{Diamondstone} generalize to infinite joins of members of effectively closed sets that do not contain computable members. 

\subsection{Lower cone avoidance}

For our purpose, we define infinite joins as follows

\begin{defn}
Let $\{A_i\}_{i\in\omega}$ be a countable sequence of sets. The {\em infinite join} is defined by
\begin{center}
$\bigoplus\{A_i\}=\{\langle i,x\rangle : x\in A_i\}$.
\end{center}
\end{defn}

\vspace{0.5cm}
We argue for that every special $\Pi^0_1$ class $\mathcal{P}$ contains a member $X$ such that $X\not\leq_T\bigoplus\{A_i\}_{i\in\omega}$ for $A_i\in\mathcal{P}$.
\vspace{0.5cm}

\begin{thm}
Let $\mathcal{P}$ be a special $\Pi^0_1$ class and let $\{A_i\}_{i\in\omega}$ be a sequence of sets in $\mathcal{P}$. Then $\mathcal{P}$ has a member that is not computable by $\bigoplus\{A_i\}_{i\in\omega}$.
\end{thm}
\begin{prf}
Let $T$ be an infinite computable tree with no computable paths, let $j\in\omega$, and let $\{A_i\}_{i\in\omega}$ be a sequence of sets in $[T]$. It suffices to show that there exists an infinite computable subtree $T^*\subset T$ such that $f\neq\Psi_j(\bigoplus\{A_i\})$ for any $f\in[T^*]$, where $f\neq A_i$.
\vspace{0.5cm}

Given that $T$ has no computable paths, $T$ contains at least two incompatible strings which are infinitely extendible. Call the least such strings $\sigma$ and $\tau$. Let $n<\textrm{min}\{|\sigma|,|\tau|\}$ be the least number such that $\sigma(n)\neq\tau(n)$. Ask the oracle for $(\bigoplus\{A_i\})'$ whether or not $\Psi_j(\bigoplus\{A_i\};n)$ converges. If so, then one of the two strings, say $\sigma$, must disagree with $\Psi_j(\bigoplus\{A_i\};n)$, and let $T^*$ consists of all strings in $T$ which are compatible with $\sigma$. We then have $f(n)=\sigma(n)$ for all $f\in [T^*]$. So $f\neq\Psi_j(\bigoplus\{A_i\})$ for all such $f$. Otherwise, we let $T^*=T$, and in this case the theorem follows trivially.\qed
\end{prf}

\subsection{Upper cone avoidance}
We now demonstrate that if $\mathcal{P}$ is a $\Pi^0_1$ class, and $X$ is a non-recursive set, then there exists a countable sequence of sets $\{A_i\}_{i\in\omega}$ in $\mathcal{P}$ such that $X\not\leq_T\bigoplus\{A_i\}_{i\in\omega}$.

\begin{thm}
Let $X$ be a non-recursive set. Every non-empty $\Pi^0_1$ class contains a countable sequence of sets whose infinite join does not compute $X$.
\end{thm}
\begin{prf}
Let $X$ be a non-recursive set, $T$ be an infinite computable tree and $i\in\omega$. We must show that there exists an infinite computable subtree $T^*\subset T$ such that $X\neq\Psi_i(\bigoplus\{A_j\})$ for any $A_j\in [T^*]$. Furthermore, an index for $T^*$ can be found $(\emptyset'\oplus X)$-recursively from $i$ and an index for $T$.
\vspace{0.5cm}

For each $n\in\omega$, define

\begin{center}
$U_n=\{\sigma\in T:\Psi_i(\oplus\sigma;n)\uparrow\vee \Psi_i(\oplus\sigma;n)\downarrow\neq X(n)\}$.
\end{center}

Note that if $U_n$ is finite, then $\oplus\sigma$ codes an element of $\omega$ since each $\sigma$ would be of finite length. If $U_n$ is infinite, then $\oplus\sigma$ codes a subset of $\omega$. We must show that some $U_n$ must have infinitely many elements.
\vspace{0.5cm}

Case 1: Every $U_n$ is finite. In this case we can find some $m\in\omega$ and some $k\in\omega$ such that $\Psi_i(\oplus\sigma;n)=k$ for all $\sigma$ such that $|\sigma |=m$. Then, it must be that $X(n)=k$ and $X$ must be recursive.
\vspace{0.5cm}

Case 2: Some $U_n$ is infinite. Using an oracle for $\emptyset'$, we find the least $n$ such that $U_n$ is infinite. We let $T^*=U_n$. As a result $\Psi_i(\bigoplus A_j)\neq X$ for every $A_j\in [T^*]$.\qed
\end{prf}

\section{Countable strict chains in the degree spectrum}

Finally, we now show that there exists a $\Pi^0_1$ class whose degree spectrum is a strictly increasing infinite sequence of Turing degrees.

\begin{defn}
Let $A=\{{\bf a}_i\}_{i\in\omega}$ be an infinite sequence of Turing degrees. We say that $A$ is a {\em strict chain} if ${\bf a}_n<{\bf a}_m\in A$ whenever $n<m$, for all ${\bf a}_n, {\bf a}_m\in A$. Let $\mathcal{P}\subseteq 2^\omega$ be an infinite $\Pi^0_1$ class. 
\end{defn}

It is known due to Jockusch and Soare \cite{JS1} that there exists a $\Pi^0_1$ class such that all members are Turing incomparable. Moreover, every special $\Pi^0_1$ class contains mutually incomparable elements. If we want to construct a class whose degree spectrum is a strict chain, our only option is to define either a countable class or, at least, a $\Pi^0_1$ class with a computable element. We shall go for the second option of constructing a $\Pi^0_1$ class containing a computable element. The rest of the elements of the class will contain only sets of r.e. degrees each of which is an element of a strictly increasing total order that can be obtained by Sacks Density Theorem \cite{SacksDensity}. Given a $\Pi^0_1$ class $\mathcal{P}$ with any degree spectrum, the reason why we are allowed to enumerate an r.e. set of degree not belonging to the degree spectrum, into $\mathcal{P}$, is due to a theorem by Cenzer and Smith \cite{CenzerSmith} which essentially says that if $\alpha$ is the degree spectrum of a $\Pi^0_1$ class then for any r.e. degree ${\bf a}$ not in $\alpha$, $\alpha\cup\{{\bf a}\}$ is the degree spectrum of a $\Pi^0_1$ class. We use a modified version of this theorem and Sacks Density Theorem to prove our claim.

\begin{thm}
There exists a $\Pi^0_1$ class $\mathcal{P}$ such that the degree spectrum of $\mathcal{P}$ is an infinite strict chain.
\end{thm}
\begin{prf}

To prove the theorem we first need to show that if $\alpha$ is the degree spectrum of a $\Pi^0_1$ class $\mathcal{P}$, then for any r.e. degree ${\bf a}\not\in\alpha$ there exists a $\Pi^0_1$ class $\mathcal{Q}$ whose degree spectrum is $\alpha\cup\{{\bf a}\}$. For this we modify the theorem by Cenzer and Smith. Although a slightly different proof was appeared in \cite{CevikChoice}, we shall modify it for the purpose of proving our result. We then apply the construction for the sequence of r.e. degrees ${\bf a}_i$ that can be obtained by a successive application of Sacks Density Theorem between the two r.e. degrees ${\bf 0}$ and ${\bf 0'}$. Applying Sacks Density will give us a countably infinite sequence of r.e. degrees $\{{\bf a}_i\}_{i\in\omega}$ such that ${\bf a}_i<{\bf a}_{i+1}$ for any $i\in\omega$. We need two definitions to proceed.

\begin{defn}
Let $A$ be an r.e. set. The {\em enumeration function} of $A$ is 
\begin{center}
$f(n)=\mu s(A_s\upharpoonright n=A\upharpoonright n)$,
\end{center}
\noindent where $A_s$ denotes the set of numbers enumerated into $A$ by the end of stage $s$.
\end{defn}

Let $f_s(n)=\mu s'\leq s$ such that $A_{s'}\upharpoonright n=A_s\upharpoonright n$. Note that $f\equiv_T A$. So if $A$ is non-computable, then so is $f$. However, $f_s$ is computable and $f_s(n)$ can only increase as $s$ increases.

\begin{defn}
A {\em copy} of a $\Pi^0_1$ class $\mathcal{P}=[\Lambda]$ for some downward closed computable set of strings $\Lambda$ is defined by the set $\{\tau*\sigma:\sigma\in\Lambda\}$ for any $\tau\in 2^{<\omega}$.
\end{defn}

To enumerate an r.e. set $A$ into a given class $\mathcal{P}$, we code the enumeration function of $A$ on a path of the new class $\mathcal{Q}$ we aim to define and put a copy of $\mathcal{P}$ above the rest of the strings. The general picture of $\mathcal{Q}$ is given in Figure 1. For the enumeration function to be coded on a path, we need to be able to code the enumeration distance between the elements that are enumerated in $A$. That is, we need to compute the number of stages between each enumerated element in $A$ as we code the enumeration function of $A$. To code the enumeration distance we put a delimeter symbol along the path of the enumeration function. Since $\mathcal{P}\subset \{0,1\}^\omega$, we should use another symbol as a delimeter not to confuse with the branches of $\mathcal{Q}$. Let $\#$ denote the delimeter symbol. Our new class $\mathcal{Q}$ will then be a subset of $\{0,1,\#\}^\omega$.

\begin{figure}[ht]
\begin{center}
\includegraphics{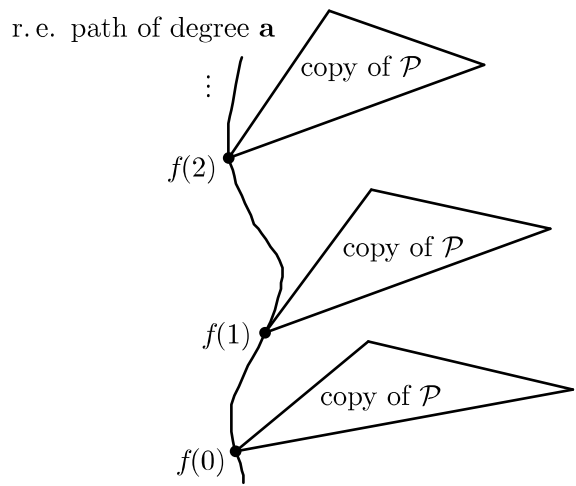}
\end{center}
\caption{$\mathcal{Q}$ contains a path of r.e. degree ${\bf a}$ and elements of degrees Turing equivalent to that of the elements of $\mathcal{P}$.}
\label{fig:figureone}
\end{figure}

Now suppose that $\mathcal{P}=[\Lambda]$ is a  $\Pi^0_1$ class with degree spectrum $\alpha$. We shall assume that we are given, initially, an infinite $\Pi^0_1$ class of computable sets. This is easy to define. At stage $s=0$, enumerate $\emptyset$ into $T_0$. At stage $s>0$, given $T_s$, for each leaf $\sigma\in T_s$ do the following: If $\sigma$ ends with 1, enumerate $\sigma*1$ into $T_s$; If $\sigma$ ends with 0, then enumerate $\sigma*0$ and $\sigma*1$ into $T_s$. Then let $T=\cup_s T_s$. Clearly $[T]$ is a countable $\Pi^0_1$ class in which all members are computable. 

\begin{figure}[ht]
\begin{center}
\includegraphics{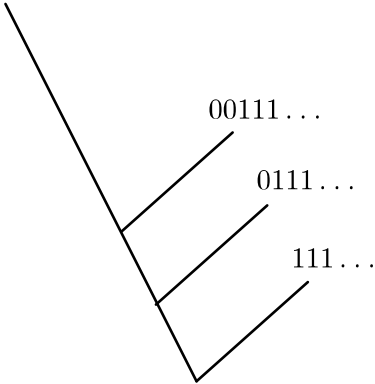}
\end{center}
\caption{A $\Pi^0_1$ class in which all members are computable.}
\label{fig:figuretwo}
\end{figure}

Suppose that we are given an r.e. degree ${\bf a}\not\in\alpha$. We define a downward closed set of strings $\Upsilon=\bigcup_{s\in\omega} \Upsilon_s$ as a subset of $\{0,1,\#\}^{<\omega}$ such that $\mathcal{Q}=[\Upsilon]$ is a $\Pi^0_1$ class with degree spectrum $\alpha\cup\{{\bf a}\}$. So $\Upsilon$ will be a ternary tree containing copies of $\mathcal{P}$ and a path of degree ${\bf a}$. For defining $\mathcal{Q}$, we start to place a copy of $\mathcal{P}$ above each delimeter symbol ($\#$ that is) in $\Upsilon_s$, in the form of a set of strings in $\Lambda$ up to length of that stage of the construction. When the value of $f_s(n)$ gets stabilized we will be leaving a whole copy of $\mathcal{P}$. We do this for every argument $n$. Therefore we will have multiple copies of $\mathcal{P}$. However, this is not a problem since we do not care whether the elements of $\mathcal{Q}$ are distinct or not. Let $\Pi_s$ be a set of strings of the form
\begin{center}
$2^{f_s(0)+1}\# 2^{f_s(1)+1}\# 2^{f_s(2)+1}\#\ldots$
\end{center}

\noindent The enumerated strings of $\Upsilon$ will be of this form.
\vspace{0.5cm}

\noindent We now give the formal construction of $\mathcal{Q}$.
\vspace{0.5cm}

{\bf Construction.}

At stage $s=0$: Enumerate $\emptyset$ into $\Upsilon_0$, define $\Upsilon^*=\emptyset$.

At each next stage $s+1$: Given $\Pi_s$ and $\Upsilon_s$, let $\sigma$ be a leaf of $\Upsilon_s$.

(i) First we enumerate $\sigma*i$ into $\Upsilon_{s+1}$ for $i\in\{0,1,\#\}$ if there exists a string $\tau\in\Pi_s$ such that $\sigma*i\subset\tau$. If $i=\#$, then enumerate $\sigma*i$ also into $\Upsilon^*$.

(ii) Secondly, to put a copy of $\mathcal{P}$, we see if there exists a string $\eta\in\Upsilon^*$ such that $\eta\subset\sigma$. If so, let $\upsilon\in\Upsilon^*$ and $\pi\in\Pi_s$ be such that $\upsilon*\pi=\sigma$. For $i\in\{0,1\}$, enumerate $\sigma*i$ into $\Upsilon_{s+1}$ if $\pi*i\in\Pi_s$.

This ends the formal construction.
\vspace{0.5cm}

\noindent Let us now verify that our construction works. We shall argue that $\mathcal{Q}$ contains a path of degree ${\bf a}$ and (possibly multiple) copies of $\mathcal{P}$. Let $s$ be an arbitrary stage and let $\sigma$ be a leaf of $\Lambda_s$. Whenever we find $\tau$ satisfying $\tau\supset\sigma*\#$, when we enumerate $\sigma*\#$ into $\Upsilon_{s+1}$ and $\Upsilon^*$, it follows that $\tau$ (without all the $\#$'s) is actually an initial segment of $\Lambda_s$. As $s$ increases, so does the enumeration distance between the elements that $f$ enumerates. It will eventually get stabilized when $f_s(n)$ is not changed anymore. Next we argue that we leave a copy of $\mathcal{P}$ in $\Upsilon$. Let $\sigma$ be the least string in $\Upsilon_s$ such that, for a sufficiently large stage $s$, $\tau\not\in\Pi_s$ for all $\tau\supset\sigma$. Such $s$ exists since every $f_s(n)$ gets changed finitely many times. Then, for all $s'\geq s$, step (ii) ensures enumerating two incompatible extensions of $\tau'$ such that $\tau*\tau'\in\Upsilon_s$ for every $\tau\in\Lambda$.
\vspace{0.5cm}

We are now ready for the final touch. Let $\{{\bf a}_i\}_{i\in\omega}$ be a countable sequence of r.e. degrees, with ${\bf a}_0$ being ${\bf 0}$, which can be obtained from a successive application of Sacks Density Theorem for r.e. degrees in between ${\bf 0}$ and ${\bf 0'}$. Apply the construction we have just given on a $\Pi^0_1$ class with all members computable and for ${\bf a}_1$ to obtain a $\Pi^0_1$ class $\mathcal{P}_1$. Then for every ${\bf a}_i$, for $i>1$ and $j\geq 1$, we inductively apply the same construction on $\mathcal{P}_j$ and ${\bf a}_{i+1}$ to obtain a $\Pi^0_1$ class $\mathcal{P}_{j+1}$. We then define $\mathcal{R}=\bigcup_{j\in\omega}\mathcal{P}_j$. Clearly, $\mathcal{R}$ is a $\Pi^0_1$ class and it satisfies the desired property.\qed
\end{prf}


\begin{thebibliography}{40}

\bibitem{CenzerRec}
D. Cenzer: $\Pi^0_1$ Classes in Recursion Theory. Handbook of Computability Theory. North-Holland, Studies in Logic and the Foundations of Mathematics Vol.140 (1999), p.37-89.




\bibitem{CenzerSmith}
D. Cenzer and R. L. Smith, On the ranked points of a $\Pi^0_1$ set. {\em Journal of Symbolic Logic}, Vol. 54, pp. 975-991 (1989).


\bibitem{Cooper}
S. B. Cooper: Computability theory, Chapman \& Hall, CRC Press, Boca Raton, FL, New York, London (2004).

\bibitem{Cevik}
A. \c{C}evik: Antibasis theorems for $\Pi^0_1$ classes and the jump hierarchy. {\em Archive for Mathematical Logic}, Vol. 52, Issue 1-2, pp. 137-142 (2013).

\bibitem{CevikChoice}
A. \c{C}evik: $\Pi^0_1$ Choice Classes. {\em Mathematical Logic Quarterly}, Vol. 62, No 6, pp. 563-574 (2016).

\bibitem{Diamondstone}
D. E. Diamondstone, D. D. Dzhafarov, R. I. Soare: $\Pi^0_1$ Classes, Peano Arithmetic, Randomness, and Computable Domination. {\em Notre Dame J. Formal Logic} Volume 51, Number 1 (2010), p.127-159. 

\bibitem{Downey}
R. Downey and D. Hirshfeldt: Algorithmic Randomness and Complexity, Springer-Verlag (2010).


\bibitem{GS97}
M. J. Groszek and T. A. Slaman: $\Pi^0_1$ classes and minimal degrees, {\em Annals of Pure and Applied Logic} {\bf 87}(2): 117-144 (1997).


\bibitem{JockuschSimpson1980}
C. G. Jockusch and S. G. Simpson: Minimal degrees, hyperimmune degrees, complete extensions of arithmetic, Abstract no. 781-E10, {\em Abstract of papers presented to the AMS} {\bf 1}: 546 (1980). 


\bibitem{JS1}
C. Jockusch, R. I. Soare, $\Pi^0_1$ classes and degrees of theories. {\em Trans. Amer. Math. Soc.}, 173, pp. 33-56 (1972).

\bibitem{JS2}
C. Jockusch, R. I. Soare: Degrees of members of $\Pi^0_1$ classes. {\em Pacific J. Math.}, 40, pp. 605-616 (1972).



\bibitem{KL}
T. Kent, A. E. M. Lewis: On the Degree Spectrum of a $\Pi^0_1$ Class. {\em Trans. Amer. Math. Soc.}, 362, pp. 5283-5319 (2010).

\bibitem{SK}
S. C. Kleene: Recursive predicates and quantifiers, {\em Trans. Amer. Math. Soc.}, Vol. 53 (1943), pp. 41-73.










\bibitem{SacksDensity}
G. E. Sacks: The recursively enumerable degrees are dense. {\em Ann. of Math.}, 80(2), pp. 300-312 (1964).




\bibitem{Soare}
R. I. Soare: Recursively Enumerable Sets and Degrees, Perspectives in Mathematical Logic. Springer-Verlag, Berlin (1987).












\end{thebibliography}
\end{document}